\documentclass [12pt]{article}
\usepackage{amssymb,amsfonts,latexsym,amsmath,amsthm}
\usepackage{comment}

\def\Z{\mathbb{Z}}

\def\strictsubset{\lower .1cm \hbox{${\buildrel{\subset}\over{_{\not=
}}}$}}

\def\resp.{{\it resp.}}

\newtheorem{theorem}{Theorem}

\newtheorem{lemma}{Lemma}

\begin{document}

\title{On special Wieferich's primes}
\author{Luis H. Gallardo}
\date{}
\maketitle

\begin{center}
Running head:  Wieferich's primes
\end{center}
\begin{center}
Mathematics, University of Brest\\
6, Avenue Le Gorgeu, C.S. 93837,\\
29238 Brest Cedex 3, France.\\
e-mail : Luis.Gallardo@univ-brest.fr\\                                                       
\end{center}
\begin{center}
\end{center}
\begin{center}
(Concerned with sequences A001220, A001348, A005384, A005385.)
\end{center}

\newpage

\begin{abstract}
We prove that there are no Wieferich's primes $q=2p+1$ where $p \equiv 3 \pmod{4}$ is a prime number.
\end{abstract}

\section{Introduction}

Let $G$ be a finite multiplicative cyclic group of cardinal $d$ with unity $1.$
Throughout the whole paper we denote by $o(G)$ the order $d$ of $G$ and we denote by
$o_d(a)$ the order of an element  $a \in G$ that is, the smallest nonnegative integer
$n \geq 0$ such that $a^n =1$ in $G.$

More precisely, given a prime number $q,$ we denote by $o_q(a)$ the order of an element $a \in (\Z/q\Z)\sp{*}$
and by $o_{q^2}(b)$ the order of an element $b \in (\Z/q^2\Z)\sp{*}.$

As usual, for positive integers $a \mid\mid b$ means that $a$ divides $b$ (noted also $a \mid b$)
and that $\gcd(a,b/a)=1.$ If $p$ is a prime number  such that $p^2 \mid 2^{p-1}-1$ then $p$ is 
{\emph{Wieferich's}} prime.  For a prime number $p$ and for a positive integer $r>0$
we denote by $(\Z/p^r\Z)\sp{*}$ the cyclic group of nonzero elements of the ring $\Z/p^r\Z.$

We call {\emph{Sophie Germain's}} prime a prime number $p$ such that $q=2p+1$ is also prime.
It is well known \cite[Theorem 103, p. 80]{hardy}
that for a Sophie Germain's prime $p \equiv 3 \pmod{4}$ the prime number $q=2p+1$ is the 
smallest prime divisor of the Mersenne number $M_p = 2^p-1.$ 
Since it is believed that Mersenne numbers $M_p$ with prime $p$ are square-free, it may have some interest to know 
whether or not $q^2$ divides $M_p.$  This can be investigated with computers since 
$(2p+1)^2$ divides $M_p$ is equivalent to $o_q(2) =p$ and this is easy to check for Sophie Germain's primes $p.$
Indeed, we checked that $q^2$ never divides $M_p$ for all Sophie Germain' primes less that $10^{11},$
in little computer time. But, of course, this was a loss of computer time, since it is well known \cite{bray}
(see also Lemma \ref{bw}) that for
prime numbers $p$
the Mersenne number $M_p$ can only have primary divisors $r^a \mid \mid M_p$ with $r$ prime and $a>1$, when 
$r$ is a Wieferich's prime, that is $2^{r-1} \equiv 1 \pmod{r^{2}}.$  Observe that for a Sophie Germain's prime
$p \equiv 3 \pmod {4}$, $q^2 \mid M_p$ implies that
$2^{2p} -1 \equiv 0 \pmod{q^2}$ so that $q=2p+1$ is a Wieferich's prime, and
Dorais and Klyve 
\cite{klyve} proved recently that
there is no Wieferich's primes less that $6.75 \cdot 10^{15}.$

Fortunately, it turns out from \cite{elehmer, lebesgue, maxfield} that it is easy to prove that $q^2$
does not divide $M_p$ when $p \equiv 3 \pmod{4}
$ is a Sophie Germain's prime. This proves also that there are no Wieferich's primes $q=2p+1$ 
where $p \equiv 3 \pmod{4}$ is a prime number. The object of this short paper is to prove these facts.

More precisely, the object of this paper is to prove the following two results:

\begin{theorem}
\label{specialw1}
Let $p \equiv 3 \pmod{4}$ be a Sophie Germain's prime. Set $q=2p+1.$ Then
\[
q \mid \mid M_p = 2^p -1.
\]
In other words, the square  $q^2$ of the smallest prime divisor $q=2p+1$ of the Mersenne number $M_p$
does not divide $M_p.$
\end{theorem}

\begin{theorem}
\label{specialw2}
Let $p \equiv 3 \pmod{4}$ be a Sophie Germain's prime. Set $q=2p+1.$ Then
$q$ is not a Wieferich's prime.
\end{theorem}

\section{Some tools}

\begin{lemma}
\label{orders}
Let $G$ be a finite cyclic group, let $g$ be a generator of $G.$
Let $x,y \in G,$ be two elements of $G.$ Let $r>0$ be a positive integer.
For an element $h \in G,$ we denote by $o(h)$ his order. We denote by $o(G)$ the order of $G.$
\begin{itemize}
\item[{\emph{a)}}]
One has $o(xy) = o(x)o(y)$ whenever $\gcd(o(x),o(y))=1.$
\item[{\emph{b)}}]
One has $o(x^r) = \frac{o(x)}{\gcd(o(x),r)}.$
\end{itemize}
\end{lemma}

Bray and Warren \cite[Theorem 1]{bray} proved that

\begin{lemma}
\label{bw}
Let $p,q$ be odd prime numbers.
\begin{itemize}
\item[{\emph{a)}}]
If $p$ divides $M_q=2^q-1$
then
\[
2^{\frac{p-1}{2}} \equiv 1 \pmod{M_q}.
\]
\item[{\emph{b)}}]
If $p^2$ divides $M_q=2^q-1$
then
\[
p^2 \mid 2^{p-1}-1.
\]
\end{itemize}
\end{lemma}

Hardy and Wright \cite[Theorem 103]{hardy} proved
and mentions Euler (see also \cite[Theorema 11, p. 28]{euler}) at the origin of:

\begin{lemma}
\label{hw}
Let $p >7$ be an odd prime number such that $p \equiv 3 \pmod{4}.$
then $q=2p+1$ is a prime number if and only if
\[
2^{p} \equiv 1 \pmod{q}.
\]
Thus, if $q$ is prime, $q \mid M_p=2^{p}-1$ and $M_p$ is composite.
\end{lemma}

Maxfield \cite[Theorem]{maxfield} proved that
\begin{lemma}
\label{mf}
Let $p$ be an odd prime number and let $r>0$ be a positive integer. 
Let $a \in G= (\Z/p\Z)\sp{*}$ be a nonzero element of $\Z/p\Z$
of order $e >1.$ Then either $a$ or $a_1$ where  $a_1 = a^{e-1}$ has order $ep^{r-1}$
in the cyclic group $G(r) =( \Z/p^r\Z)\sp{*}.$
\end{lemma}

The special case when $a$ generates $G$ was first announced by Lebesgue \cite{lebesgue}:

\begin{lemma}
\label{vl}
Let $p$ be an odd prime number and let $r>0$ be a positive integer. 
Let $a \in G= (\Z/p\Z)\sp{*}$ be a generator of $G$, that is $o_p(a)=p-1.$
Then either $a$ or $a_1$ where  $a_1 = a^{p-2}$ generates
the cyclic group $G(r) =( \Z/p^r\Z)\sp{*}.$
\end{lemma}

As usual define the Bernoulli numbers $B_{\nu}$ by:
\begin{gather}
B_0 = 1, B_1 = -\frac{1}{2}, B_2 = \frac{1}{6}, B_{2k+1}=0 \quad \text{for} \quad k>0,\\
B_r = \sum_{n=0}^{r} {r \choose n} B_n
\end{gather}

\noindent Emma Lehmer \cite{elehmer} proved that

\begin{lemma}
\label{el}
Let $p$ be an odd prime number. Let $\nu >0$ be a positive integer such that
$\nu \not \equiv 1 \pmod{p - 1}.$ Then
\[
\sum_{r=1}^{p-1} r^{\nu} \equiv pB_{\nu} \pmod{p^2}
\]
\end{lemma}

\section{Proof of the Theorems}
\label{mersenne}

Theorem \ref{specialw2} follows from Theorem \ref{specialw1} since by Lemma \ref{hw} we have $q \nmid 2^p+1.$
In more detail: Assume that $q=2p+1$ is a Wieferich's prime.  Then $q^2 \mid (2^p+1)(2^p-1)$
but by Lemma \ref{hw},  $q \nmid 2^p+1,$ so that $q^2 \mid 2^p-1.$ But this contradicts Theorem \ref{specialw1}.

We prove now Theorem \ref{specialw1}.  Observe that $p \geq 11,$ so that $q>3.$ By Lemma \ref{hw} we obtain $q \mid M_p.$
Assume now that $q^2 \mid M_p.$ Thus, $2^q \equiv 2 \pmod{q^2}$ and $o_{q^2}(2)=p.$ We shall produce a contradiction.

 First proof is as follows:
Observe that $x= -2$ has order 
\[
o_{q}(x) = o_{q}(-1)o_{q}(2) =2o_{q}(2)
\]
 by Lemma \ref{orders}. But $o_{q}(2)=p,$ so $o_q(x)=2p.$
In other words $x$ generate the cyclic group $(\Z/q\Z)\sp{*}.$ So $\{x,x^2,\ldots,x^{q-1}\} = \{1,2,\ldots,q-1\}.$ Thus computing
both sides of the congruence modulo $q^2$  of Lemma \ref{el} with $r=x$  and $\nu =2$ we get the contradiction:
\[
1 \equiv  \frac{4^q -1}{3} \equiv \sum_{k=0}^{q-1} {((-2)^k)^2} \equiv \sum_{r=1}^{q-1} r^2 \equiv qB_2 \equiv \frac{q}{6} \pmod{q^2}.
\]

A second proof is as follows:  Observe that 
\[r=o((\Z/q^2\Z)\sp{*})=2pq.
\]
As before, $x=-2$ generates $(\Z/q\Z)\sp{*}.$ But $x^{2p} \equiv 2^{2p} \equiv 1 \pmod{q^2}$
so by Lemma \ref{vl}, $y \equiv x^{q-2} \pmod{q^2}$ must generate $(\Z/q^2\Z)\sp{*}.$ But this is impossible since
$y \equiv x^{q-2} \equiv  -\frac{2^q}{4} \equiv -\frac{1}{2} \pmod{q^2}$ and
$y^{2p} \equiv 1 \pmod{q^2}.$ Alternatively, from Lemma \ref{orders} b) we get the contradiction $o_{q^2}(y)=2p$
since $o_{q^2}(-2) =2p$ and $q-2=2p-1.$

A third proof comes from Lemma \ref{mf}: Observe that $p=o_q(2).$ So by Lemma \ref{mf} either $x=2$ or $y=2^{p-1}$
has order $pq$ in the cyclic group $(\Z/q^2\Z)\sp{*}.$ But both statements are false since $p = o_{q^2}(2).$

\def\thebibliography#1{\section*{\titrebibliographie}
\addcontentsline{toc}
{section}{\titrebibliographie}\list{[\arabic{enumi}]}{\settowidth
 \labelwidth{[
#1]}\leftmargin\labelwidth \advance\leftmargin\labelsep
\usecounter{enumi}}
\def\newblock{\hskip .11em plus .33em minus -.07em} \sloppy
\sfcode`\.=1000\relax}
\let\endthebibliography=\endlist
\def\biblio{\def\titrebibliographie{References}\thebibliography}
\let\endbiblio=\endthebibliography

%%%% MACROS DE SEROUL POUR LES REFERENCES %%%%

%%%%%%% bibliographie selon AMS style %%%%%%%%%%%
%%%%%%% inspir de TUGboat 11 (1990), p. 609 %%%%%%%

\newbox\auteurbox
\newbox\titrebox
\newbox\titrelbox
\newbox\editeurbox
\newbox\anneebox
\newbox\anneelbox
\newbox\journalbox
\newbox\volumebox
\newbox\pagesbox
\newbox\diversbox
\newbox\collectionbox
%--------------------------------------------
\def\fabriquebox#1#2{\par\egroup
\setbox#1=\vbox\bgroup \leftskip=0pt \hsize=\maxdimen \noindent#2}
%--------------------------------------------
\def\bibref#1{\bibitem{#1}

%\mbox{}\ignorespaces
% warn: do not uncomment line above

\setbox0=\vbox\bgroup}
%--------------------------------------------
\def\auteur{\fabriquebox\auteurbox\styleauteur}
\def\titre{\fabriquebox\titrebox\styletitre}
\def\titrelivre{\fabriquebox\titrelbox\styletitrelivre}
\def\editeur{\fabriquebox\editeurbox\styleediteur}

\def\journal{\fabriquebox\journalbox\stylejournal}

\def\volume{\fabriquebox\volumebox\stylevolume}
\def\collection{\fabriquebox\collectionbox\stylecollection}
%--------------------------------------------
{\catcode`\- =\active\gdef\annee{\fabriquebox\anneebox\catcode`\-
=\active\def -{\hbox{\rm
\string-\string-}}\styleannee\ignorespaces}}
%--------------------------------------------
{\catcode`\-
=\active\gdef\anneelivre{\fabriquebox\anneelbox\catcode`\-=
\active\def-{\hbox{\rm \string-\string-}}\styleanneelivre}}
%--------------------------------------------
{\catcode`\-=\active\gdef\pages{\fabriquebox\pagesbox\catcode`\-
=\active\def -{\hbox{\rm\string-\string-}}\stylepages}}
%--------------------------------------------
{\catcode`\-
=\active\gdef\divers{\fabriquebox\diversbox\catcode`\-=\active
\def-{\hbox{\rm\string-\string-}}\rm}}
%--------------------------------------------
\def\ajoutref#1{\setbox0=\vbox{\unvbox#1\global\setbox1=
\lastbox}\unhbox1 \unskip\unskip\unpenalty}
%--------------------------------------------
\newif\ifpreviousitem
\global\previousitemfalse
\def\separateur{\ifpreviousitem {,\ }\fi}
%--------------------------------------------
\def\voidallboxes
{\setbox0=\box\auteurbox \setbox0=\box\titrebox
\setbox0=\box\titrelbox \setbox0=\box\editeurbox
\setbox0=\box\anneebox \setbox0=\box\anneelbox
\setbox0=\box\journalbox \setbox0=\box\volumebox
\setbox0=\box\pagesbox \setbox0=\box\diversbox
\setbox0=\box\collectionbox \setbox0=\null}
%--------------------------------------------
\def\fabriquelivre
{\ifdim\ht\auteurbox>0pt
\ajoutref\auteurbox\global\previousitemtrue\fi
\ifdim\ht\titrelbox>0pt
\separateur\ajoutref\titrelbox\global\previousitemtrue\fi
\ifdim\ht\collectionbox>0pt
\separateur\ajoutref\collectionbox\global\previousitemtrue\fi
\ifdim\ht\editeurbox>0pt
\separateur\ajoutref\editeurbox\global\previousitemtrue\fi
\ifdim\ht\anneelbox>0pt \separateur \ajoutref\anneelbox
\fi\global\previousitemfalse}
%--------------------------------------------
\def\fabriquearticle
{\ifdim\ht\auteurbox>0pt        \ajoutref\auteurbox
\global\previousitemtrue\fi \ifdim\ht\titrebox>0pt
\separateur\ajoutref\titrebox\global\previousitemtrue\fi
\ifdim\ht\titrelbox>0pt \separateur{\rm in}\
\ajoutref\titrelbox\global \previousitemtrue\fi
\ifdim\ht\journalbox>0pt \separateur
\ajoutref\journalbox\global\previousitemtrue\fi
\ifdim\ht\volumebox>0pt \ \ajoutref\volumebox\fi
\ifdim\ht\anneebox>0pt  \ {\rm(}\ajoutref\anneebox \rm)\fi
\ifdim\ht\pagesbox>0pt
\separateur\ajoutref\pagesbox\fi\global\previousitemfalse}
%--------------------------------------------
\def\fabriquedivers
{\ifdim\ht\auteurbox>0pt
\ajoutref\auteurbox\global\previousitemtrue\fi
\ifdim\ht\diversbox>0pt \separateur\ajoutref\diversbox\fi}
%--------------------------------------------
\def\endbibref
{\egroup \ifdim\ht\journalbox>0pt \fabriquearticle
\else\ifdim\ht\editeurbox>0pt \fabriquelivre
\else\ifdim\ht\diversbox>0pt \fabriquedivers \fi\fi\fi
.\voidallboxes}
%--------------------------------------------

\let\styleauteur=\sc
\let\styletitre=\it
\let\styletitrelivre=\sl
\let\stylejournal=\rm
\let\stylevolume=\bf
\let\styleannee=\rm
\let\stylepages=\rm
\let\stylecollection=\rm
\let\styleediteur=\rm
\let\styleanneelivre=\rm

\begin{biblio}{99}

\begin{bibref}{bray}
\auteur{H. G. Bray, L. R. Warren} 
\titre{On the square-freeness of Fermat and Mersenne numbers}
\journal{Pacific J. Math.} \volume{22} \pages 563-564 \annee 1967
\end{bibref}

\begin{bibref}{klyve}
\auteur{F. G. Dorais, D. W. Klyve} 
\titre{Near Wieferich Primes up to $6.7 \times 10^{15}$}
\journal{Available online at http://www-personal.umich.edu/{} \textasciitilde dorais/docs/wieferich.pdf }
\volume{} \pages  \annee 2008
\end{bibref}

\begin{bibref}{euler}
\auteur{L.  Euler}
\titre{Theoremata circa divisores numerorum, Novi Commentarii academiae
scientiarum Petropolitanae,
Reprinted in}\journal{ Euler archive  [E134].
http://www.math.dartmouth.edu/{} \textasciitilde euler}
\volume{1} \pages 20-48 \annee 1750
\end{bibref}

\begin{bibref}{hardy}
\auteur{G. H. Hardy, E. M. Wright}
\titrelivre{An introduction to the theory of numbers - 4th Edit.}
\editeur{Oxford : At the Clarendon Press - XVI - 421 p.}
\anneelivre 1960
\end{bibref}

\begin{bibref}{lebesgue}
\auteur{V. A. Lebesgue} \titre{Th\'eor\`eme sur les racines primitives}
\journal{Comptes Rendus} \volume{64} \pages 1268-1269 \annee 1867
\end{bibref}

\begin{bibref}{elehmer}
\auteur{E. Lehmer} \titre{On congruences involving Bernoulli numbers and the quotients of Fermat and Wilson}
\journal{Ann. of Math.} \volume{39} \pages 350-360 \annee 1938
\end{bibref}

\begin{bibref}{maxfield}
\auteur{J. Maxfield, M. Maxfield} \titre{The existence of integers less than $p$ belonging to $ep^{r-1} \pmod{p}$}
\journal{Math. Mag.} \volume{33} \pages 219-220 \annee 1959/1960
\end{bibref}

\end{biblio}

\end{document}